\documentclass[11pt]{amsart}
\usepackage{amsmath, natbib}
\numberwithin{equation}{section}

\newtheorem{lemma}{Lemma}[section]

\newtheorem{theorem}[lemma]{Theorem}

\newtheorem{corollary}[lemma]{Corollary}
\newtheorem{example}[lemma]{Example}

\newtheorem{remark}[lemma]{Remark}

\newcommand{\SaS}{S$\alpha$S}

\newcommand{\BA}{{\bf A}}
\newcommand{\BB}{{\bf B}}
\newcommand{\BX}{{\bf X}}
\newcommand{\BY}{{\bf Y}}

\newcommand{\bbb}{{\mathcal B}}
\newcommand{\bbw}{{\mathcal W}}

\newcommand{\one}{{\bf 1}}
\newcommand{\eid}{\buildrel{\rm d}\over {=}}

\newcommand{\field}[1]{\mathbb{#1}}

\newcommand{\ppw}{{\tt P}_w}
\newcommand{\nnw}{{\tt N}_w}

\newcommand{\reals}{{I\!\!R}}
\newcommand{\bbr}{\reals}

\begin{document}

\title[Structure of stationary stable processes]{Null flows, positive
flows and the structure of stationary symmetric stable processes}

\author[G. Samorodnitsky]{Gennady Samorodnitsky}
\address{School of Operations Research and Industrial Engineering\\
Cornell University \\
Ithaca, NY 14853}
\email{gennady@orie.cornell.edu}
\urladdr{www.orie.cornell.edu/${\,}\tilde{}$gennady}

\thanks{\textup{2000} {\it Mathematics Subject
     Classification}, Primary  {60G10, 60G52}; Secondary 37A40.}
\thanks{Research partially supported by NSF grant DMS-0303493 and NSA
grant MSPF-02G-183  at Cornell University.}


\begin{abstract}
\noindent
This paper elucidates the connection between stationary symmetric
$\alpha$-stable processes with $0<\alpha<2$ and nonsingular flows on
measure spaces by describing a new and unique decomposition of
stationary stable processes into those corresponding 
to positive flows and those corresponding to null flows. We show that a
necessary and sufficient for a stationary stable process to be ergodic
is that its positive component vanishes. 
\end{abstract}

\maketitle

\baselineskip=18pt

\section{Introduction}\label{s:intro}

We consider the class of measurable stationary symmetric
$\alpha$-stable (henceforth, \SaS) processes $\BX=(X(t),\, t\in T)$,
$0<\alpha<2$, both in discrete time ($T=\field{Z}$) and in continuous
time ($T=\field{R}$). Unlike their Gaussian counterparts ($\alpha=2$)
which can be neatly described by either covariance function or by the
spectral measure, the structure of \SaS\ processes with $0<\alpha<2$
remained largely a mystery until the paper of \cite{rosinski:1995}
which showed that integral representations of stationary symmetric
\SaS\ processes can be taken to be of a special, and particularly
illuminating, form. The general integral representations of
$\alpha$-stable processes, of the type 
\begin{equation} \label{e:stab.process}
X(t) = \int_E f_t(x) \; M(dx), \; t\in T\, 
\end{equation}
where $M$ is a \SaS\ random measure on $E$ with a $\sigma$-finite
control measure $m$, and $f_t\in L^\alpha(m)$ for each $t$ have been
known since \cite{bretagnolle:dacunha-castelle:krivine:1966} and
\cite{schreiber:1972}; see also \cite{schilder:1970} and
\cite{kuelbs:1973} (for basic information on integrals with respect to
stable random measures one can consult
\cite{samorodnitsky:taqqu:1994}). 
What \cite{rosinski:1995} showed was
that for  measurable stationary \SaS\ processes one can
choose the kernel $f_t, t\in T$ in (\ref{e:stab.process}) in a special
way, described below. Note that we are considering both real-valued
and complex-valued \SaS\ processes. In the complex-valued case we are
assuming that the stable process is {\it isotropic}, or {\it
rotationally invariant} (see \cite{samorodnitsky:taqqu:1994}). In
that case the random measure $M$ in the integral representation can
chosen to be isotropic as well.  

According to \cite{rosinski:1995} one can choose the kernel in
(\ref{e:stab.process}) in the form 
\begin{equation} \label{e:station.kernel}
f_t(x) = a_t(x) \,\left( \frac{dm\circ
\phi_t}{dm}(x)\right)^{1/\alpha}  f\circ \phi_t(x)
\end{equation}
for $t\in T$ and $x\in E$, where $f$ is a {\it single} function in
$L^\alpha(m)$. Here $(\phi_t)$ is a measurable family of maps from $E$
onto $E$ such that 
$\phi_{t+s}(x) =\phi_t(\phi_s(x))$ for all $t,s\in\reals$ and $x\in
E$, $\phi_0(x)=x$ for all  $x\in E$, and $m\circ \phi_t^{-1}\sim m$
for all $t\in \reals$. The assumptions mean that the family $(\phi_t)$
forms a measurable nonsingular flow on $E$. Finally, $(a_t)$ is a
measurable family of $\{-1,1\}$-valued (unit circle-valued in the
complex case) functions on $E$ such that for
every $s,t\in\reals$ we have $a_{t+s}(x) = a_s(x)a_t(\phi_s(x))$
$m$-almost everywhere on $E$. This means that the family $(a_t)$ forms
a cocycle for the flow $(\phi_t)$.

One can, therefore, relate, in principle,  the ergodic-theoretical
properties of the 
flow $(\phi_t)$ to the probabilistic properties of the \SaS\ process 
$(X(t),\, t\in T)$, and under the lack of redundancy assumption 
\begin{equation} \label{e:full.support}
\text{supp}\{ f\circ \phi_t,\, t\in T\} = E,
\end{equation}
(which can, obviously, be always assumed) certain important
ergodic-theoretical properties of the flow $(\phi_t)$ remain permanent
for a given process. That is, they do not change from a representation
to representation. 

One such permanent feature of stationary \SaS\ processes is the
{\it conservative-dissipative} decomposition. Recall the Hopf
decomposition of the flow $E=C\cup D$, where $C$ and $D$ are disjoint
flow invariant measurable sets, such that the flow is conservative on
$C$ and dissipative on $D$ (see e.g. \cite{krengel:1985}). Then by
writing 
\begin{equation} \label{e:diss-cons}
X(t) = \int_C f_t(x) \; M(dx) + \int_D f_t(x) \; M(dx)
:= X^C(t) + X^D(t), \; t\in T\, 
\end{equation}
one obtains a unique in law decomposition of a stationary \SaS\
process into a sum of two independent stationary \SaS\ processes, one
corresponding to a conservative flow and the to a dissipative flow
(\cite{rosinski:1995}). (As far as terminology is concerned, we will
say that a process corresponds to a particular flow, or is generated
by a particular flow, if the process has an integral representation
(\ref{e:stab.process}) with the kernel of the form (\ref{e:station.kernel})
and (\ref{e:full.support}) holds.)

Since the processes corresponding to a dissipative flow have an
intuitively clear representation as mixed moving averages (introduced
by \cite{surgailis:rosinski:mandrekar:cambanis:1993}), attention
focused on discerning further structures in the class of stationary
\SaS\ processes generated by conservative flows. The class of {\it
cyclic} processes that generalize harmonizable processes was
introduced and discussed by \cite{pipiras:taqqu:2004}. Other examples
were discussed in \cite{rosinski:samorodnitsky:1996} and
\cite{samorodnitsky:2004c}. 

A further challenge has been to find explicit relations between the
ergodic-theoretical properties of the flow underlying a stationary
\SaS\ process and probabilistic properties of the process. Certain
steps in this direction has been taken, for example, in
\cite{resnick:samorodnitsky:xue:1999},
\cite{resnick:samorodnitsky:xue:2000} and
\cite{mikosch:samorodnitsky:2000a}. The most explicit connection of
this kind established up to date is the fact that the partial maxima
of the processes generated by conservative flows grow strictly slower
than those of the processes generated by dissipative flows; see
\cite{samorodnitsky:2004a} and \cite{samorodnitsky:2004c}. 

In this paper we present a decomposition of stationary
\SaS\ processes different from the conservative-dissipative
decomposition (\ref{e:diss-cons}), the {\it positive-null}
decomposition. While process generated by a dissipative flow is
generated by a null flow, processes generated by conservative flows
can be generated by either positive or null flow. This will,
therefore, clarify further the structure of stationary
\SaS\ processes generated by conservative flows. We will see, further,
that a stationary \SaS\ process is ergodic if and only if it is
generated by a null flow. 

We finish this introductory section by mentioning that the assumption
of measurability of the processes we are considering is automatic in
the discrete time case, and is equivalent 
to the assumption of continuity in probability in the continuous time
case. See Section 1.6 \cite{aaronson:1997} or
\cite{pipiras:taqqu:2004}.  

The following section provides the basic decomposition of a stationary
\SaS\ process. Section \ref{s:ergodicity} relates this decomposition
to ergodicity of the stable process. Finally, Section \ref{s:examples}
discusses a number of examples of stationary stable processes in the
context of the theory developed in this paper.

\section{Decomposition of stationary \SaS\ processes: positive and
null parts} \label{s:null-pos}

A nonsingular map $\phi$ on a $\sigma$-finite measure space $(E, {\mathcal
E}, m)$ is a one-to-one map $E\to E$ such that both $\phi$ and
$\phi^{-1}$ are measurable and map the measure $m$ into equivalent
measures. A nonsingular map $\phi$ is called {\it positive} if there is a
finite measure $\mu$ on $(E, {\mathcal E})$ equivalent to $m$
that is preserved under $\phi$. A subset $B$ of $E$ is called weakly
wandering if there is a sequence $0=n_0<n_1<n_2<\ldots $ such that the
sets $\phi^{-n_k} B, \, k=0,1,2, \ldots $ are disjoint. The positive-null
decomposition of the map $\phi$ is a decomposition of $E$ into a union
of two disjoint $\phi$-invariant measurable sets $\tt P$ and $\tt N$
(the positive and null parts of $\phi$) 
such that the restriction of $\phi$ to $\tt P$ is positive, and 
${\tt N} = \cup_{k=0}^\infty \phi^{-n_k} B$, with $B$ a measurable 
weakly wandering set and a  sequence $0=n_0<n_1<n_2<\ldots $ as above.
This decomposition is unique modulo null sets. See Section 1.4 in 
\cite{aaronson:1997} and Section 3.4 in \cite{krengel:1985}. (In
general, the above two references are a source of the basic facts in
ergodic theory used in this paper.)

Given a (measurable) nonsingular flow $(\phi_t)$, each map $\phi_t$ is
a nonsingular map, and has a corresponding positive-null decomposition
$E={\tt P}_t\cup {\tt N}_t$. One can show (\cite{krengel:1985}) that
there are invariant measurable sets $\tt P$ and $\tt N$ (the positive
and null parts of the flow) such that
${\tt P}_t={\tt P}$ and ${\tt N}_t={\tt N}$ modulo null sets for all
$t>0$. The flow is positive if ${\tt P}=E$ and null if ${\tt N}=E$
modulo null sets. 

Our first result shows the invariant nature of the property of being
generated by a positive or by a null flow, and is parallel to the
corresponding result for conservative and dissipative flows in
\cite{rosinski:1995}. Let $\bbw$ be class of functions $w:T\to
[0,\infty)$ such that $w$ is nonincreasing on $T\cap (-\infty,0]$,
nondecreasing on $T\cap [0,\infty)$ and 
\begin{equation} \label{e:weights}
\int_{T\cap (-\infty,0]} w(t)\lambda(dt) = \int_{T\cap [0,\infty)}
  w(t)\lambda(dt) =\infty.
\end{equation} 
Here $\lambda$ is the counting measure is $T=\field{Z}$ and the
Lebesgue measure if $T=\field{R}$. The (implicit) assumption in the
theorem below is that the process is given in the form
(\ref{e:stab.process}) with the kernel of the form
(\ref{e:station.kernel}) and (\ref{e:full.support}) holds.
\begin{theorem} \label{t:invariant}
(i) A stationary \SaS\ process is generated by a positive flow if and
 only for any $w\in \bbw$ we have
\begin{equation} \label{e:w.int.infinite}
\int_T w(t) |f|^\alpha\circ \phi_t \frac{dm\circ
\phi_t}{dm} \lambda(dt)=\infty \ \ \text{$m$-a.e.}
\end{equation}

(ii) A stationary \SaS\ process is generated by a null flow if and
 only for some $w\in \bbw$ we have
\begin{equation} \label{e:w.int.finite}
\int_T w(t) |f|^\alpha\circ \phi_t \frac{dm\circ
\phi_t}{dm} \lambda(dt)<\infty \ \ \text{$m$-a.e.}
\end{equation}

(iii) If a stationary \SaS\ process is generated by a positive
(respectively, null) flow in one representation, then in every other
representation it will 
be generated by a positive (respectively, null) flow. In particular,
the classes of stationary \SaS\ processes generated by positive and
null flows are disjoint.
\end{theorem}
\begin{proof}
(i) We start with the discrete time case $T=\field{Z}$. Suppose that
  (\ref{e:w.int.infinite}) holds for all $w\in \bbw$ and assume that,
  to the contrary, $\tt N$ is not empty $m$-a.e. For a nonsingular map
  $\phi$ and a 
nonincreasing nonnegative nonsummable sequence $w=(w_n, \,
n=0,1,2,\ldots )$  there is a (uniquely determined modulo null sets)
  invariant set $\ppw^\phi$ such that
for every strictly positive $g\in L^1(m)$ 
\begin{equation} \label{e:pos.w}
\sum_{n=0}^\infty w(n) |g\circ \phi_n(x)|^\alpha \frac{dm\circ
\phi_n}{dm}(x) =\infty \ \ \text{$m$-a.e. on $\ppw^\phi$}
\end{equation}
and 
\begin{equation} \label{e:null.w}
\sum_{n=0}^\infty w(n) |g\circ \phi_n(x)|^\alpha \frac{dm\circ
\phi_n}{dm}(x) <\infty \ \ \text{$m$-a.e. on $\nnw^\phi=
E-\ppw^\phi$,} 
\end{equation}
see Theorem 2 in \cite{krengel:1967}. 

It follows by Theorems 2 and 3 of \cite{krengel:1967} that there is a
nonincreasing nonnegative nonsummable sequence $w^{(+)}=(w_n^{(+)}, \,
n=0,1,2,\ldots )$  such that ${\tt N}_{w^{(+)}}^\phi= {\tt N}$ $m$-a.e. 

Note, furthermore, that the null parts of $\phi$ and $\phi^{-1}$
coincide modulo null sets. Therefore, appealing once again to
\cite{krengel:1967}  we see that there is a
nonincreasing nonnegative nonsummable sequence $w^{(-)}=(w_n^{(-)}, \,
n=0,1,2,\ldots )$  such that ${\tt N}_{w^{(-)}}^{\phi^{-1}}= {\tt N}$
$m$-a.e.  

Defining $w_0=w_0^{(+)}+ w_0^{(-)}$, $w_n = w_n^{(+)}$ for $n\geq 1$
and $w_n = w_n^{(-)}$ for $n\leq -1$ we obtain a $w\in \bbw$ for which
(\ref{e:w.int.infinite}) fails at almost every point of $\tt N$. Hence
$\phi$ is positive. 

In the opposite direction, suppose that $\phi$ is positive. By 
Theorem 2  of \cite{krengel:1967}, for every $w\in \bbw$,
$$
\sum_{n=0}^\infty w(n) |f\circ \phi_n|^\alpha \frac{dm\circ
\phi_n}{dm}=\infty \ \ \text{$m$-a.e. on} \ \ 
\left\{ \sum_{n=0}^\infty  |f\circ \phi_n|^\alpha >0\right\}
$$
and 
$$
\sum_{n=-\infty}^0 w(n) |f\circ \phi_n|^\alpha \frac{dm\circ
\phi_n}{dm}=\infty \ \ \text{$m$-a.e. on} \ \ 
\left\{ \sum_{n=-\infty}^0 |f\circ \phi_n|^\alpha >0\right\}.
$$
Since 
$$
E = \left\{ \sum_{n=0}^\infty  |f\circ \phi_n|^\alpha >0\right\} \cup
\left\{ \sum_{n=0}^\infty  |f\circ \phi_n|^\alpha >0\right\}
$$
$m$-a.e., we obtain (\ref{e:w.int.infinite}). This finishes the proof
of part (i) in the case $T=\field{Z}$. 

Consider now the case $T=\field{R}$. Suppose  that the flow $(\phi_t)$
is positive, and let $w\in\bbw$. We have 
\begin{equation} \label{e:cts.disc1}
\int_0^\infty w(t) |f|^\alpha\circ \phi_t \frac{dm\circ
\phi_t}{dm} \lambda(dt)
\geq \sum_{j=1}^\infty w(j)\int_0^1 |f|^\alpha\circ \phi_{t+j-1}
\frac{dm\circ
\phi_{t+j-1}}{dm} \lambda(dt)
\end{equation}
$$
= \sum_{j=1}^\infty w(j) g_+\circ \phi_j \frac{dm\circ
\phi_j}{dm},
$$
where 
\begin{equation} \label{e:gplus}
g_+(x) = \int_0^1 |f\circ \phi_{t-1}(x)|^\alpha
\frac{dm\circ
\phi_{t-1}}{dm}(x) \lambda(dt), \ x\in E.
\end{equation}
Therefore, by Theorem 2  of \cite{krengel:1967},
\begin{equation} \label{e:cts.disc2}
\int_0^\infty w(t) |f|^\alpha\circ \phi_t \frac{dm\circ
\phi_t}{dm} \lambda(dt) = \infty \ \ \text{$m$-a.e. on}
\ \ \left\{ \sum_{j=0}^\infty g_+\circ \phi_j >0\right\}.
\end{equation}
Similarly,
\begin{equation} \label{e:cts.disc3}
\int_{-\infty}^0 w(t) |f|^\alpha\circ \phi_t \frac{dm\circ
\phi_t}{dm} \lambda(dt) = \infty \ \ \text{$m$-a.e. on}
\ \ \left\{ \sum_{j=0}^\infty g_-\circ \phi_{-j} >0\right\},
\end{equation}
where 
\begin{equation} \label{e:gminus}
g_-(x) = \int_0^1 |f\circ \phi_{-t+1}(x)|^\alpha
\frac{dm\circ
\phi_{-t+1}}{dm}(x) \lambda(dt), \ x\in E.
\end{equation}
Therefore, to prove (\ref{e:w.int.infinite}) we need to show that
\begin{equation} \label{e:cts.disc4}
\left\{ \sum_{j=0}^\infty g_+\circ \phi_j >0\right\}
\cup \left\{ \sum_{j=0}^\infty g_-\circ \phi_{-j} >0\right\}
= E \ \ \text{$m$-a.e.}
\end{equation}
To this end, note that
$$
\left\{ \sum_{j=0}^\infty g_+\circ \phi_j >0\right\}
\cup \left\{ \sum_{j=0}^\infty g_-\circ \phi_{-j} >0\right\}
$$
$$
= \left\{ \sum_{j=0}^\infty g_+\circ \phi_j \frac{dm\circ
\phi_j}{dm} >0\right\}
\cup \left\{ \sum_{j=0}^\infty g_-\circ \phi_{-j} \frac{dm\circ
\phi_{-j}}{dm} >0\right\} \ \ \text{$m$-a.e.}
$$
and that 
$$
\sum_{j=0}^\infty g_+\circ \phi_j \frac{dm\circ
\phi_j}{dm}
= \int_{-1}^\infty |f|^\alpha\circ \phi_t \frac{dm\circ
\phi_t}{dm} \lambda(dt)
$$
and
$$
\sum_{j=0}^\infty g_-\circ \phi_{-j} \frac{dm\circ
\phi_{-j}}{dm}
= \int^{1}_{-\infty} |f|^\alpha\circ \phi_t \frac{dm\circ
\phi_t}{dm} \lambda(dt). 
$$
If (\ref{e:cts.disc4}) fails, then for some set $A$ with $m(A)>0$ we
have, for all $x\in A$, $f_t(x)=0$ for almost every $t\in \field{R}$.
By Fubini's theorem we conclude that for almost every $t\in \field{R}$
we have $f_t(x)=0$ $m$-a.e. on $A$. Let now $t\in \field{R}$, and
choose a sequence $t_n\to t$ such that $f_{t_n}(x)=0$ $m$-a.e. on $A$
for each $n$. By the continuity in probability
$$
0 = \lim_{n\to\infty} \int_E\left| f_t(x)-f_{t_n}(x)\right|^\alpha
m(dx)
\geq \int_E\left| f_t(x)\right|^\alpha m(dx),
$$
and for each $t\in \field{R}$  we have $f_t(x)=0$ $m$-a.e. on $A$,
contradicting the assumption of full support (\ref{e:full.support}). 
Therefore, (\ref{e:w.int.infinite}) follows. 

In the opposite direction, suppose that (\ref{e:w.int.infinite}) holds
  for all $w\in \bbw$ and assume that,   to the contrary, $\tt N$ is
  not empty $m$-a.e. Applying, once again, to Theorems 2 and 3 of
  \cite{krengel:1967} we conclude that there a nonincreasing
  nonnegative nonsummable sequence $w^{(+)}=(w_n^{(+)}, \, 
n=0,1,2,\ldots )$  such that
$$
\sum_{n=0}^\infty w(n) |g_+\circ \phi_n(x)|^\alpha \frac{dm\circ
\phi_n}{dm}(x) <\infty \ \ \text{a.e. on {\tt N}}
$$
with the function $g_+$ given in (\ref{e:gplus}). Letting $w(t) =
w_n^{(+)}$ for $n-1\leq t<n$, $n=2,3,\ldots $ and $w(t) =
w_0^{(+)}$ for $0<t<1$, we see that 
\begin{equation} \label{e:cts.disc5}
\int_0^\infty w(t) |f|^\alpha\circ \phi_t \frac{dm\circ
\phi_t}{dm} \lambda(dt) <\infty \ \ \text{a.e. on {\tt N}.}
\end{equation}
Similarly we can find a nonincreasing
  nonnegative nonsummable sequence $w^{(-)}=(w_n^{(-)}, \, 
n=0,1,2,\ldots )$  such that with $w(-t) =
w_n^{(-)}$ for $n-1\leq t<n$, $n=2,3,\ldots $ and $w(-t) =
w_0^{(-)}$ for $0<t<1$, we have
\begin{equation} \label{e:cts.disc6}
\int_{-\infty}^0 w(t) |f|^\alpha\circ \phi_t \frac{dm\circ
\phi_t}{dm} \lambda(dt) < \infty \ \ \text{a.e. on {\tt N}.}
\end{equation}
If we select $w(0)$ large enough, we obtain $w\in \bbw$ for which 
(\ref{e:w.int.infinite}) fails at a.e. point of $\tt N$. This
contradiction finishes the proof of part (i) in all cases. 

(ii) Consider first the case $T=\field{Z}$. If (\ref{e:w.int.finite})
holds for some $w\in \bbw$, then the assumption of full support
(\ref{e:full.support}) and Theorem 2 in \cite{krengel:1967} show that
${\tt P}=\emptyset$ and so $\phi$ is null. If, on the other hand,
$\phi$ is null, then appealing, as above, to Theorem 3 in
\cite{krengel:1967} we can construct $w\in \bbw$ for which
(\ref{e:w.int.finite}) holds . In the case $T=\field{R}$ we reduce the
argument to the discrete time case as in the proof of (i).

(iii) Suppose that $\BX$ is generated by a positive flow in the
representation (\ref{e:stab.process}) with the kernel of the form
(\ref{e:station.kernel}), and consider another representation (in law)
of $\BX$
$$
X(t) = \int_S g_t(y) \; M_1(dy), \; t\in T,
$$
where $M_1$ has a control measure $m_1$ and 
$$
g_t(y) = a^{(1)}_t(y) \,\left( \frac{dm_1\circ
\psi_t}{dm_1}(y)\right)^{1/\alpha}  g\circ \psi_t(y)
$$
with $(\psi_t)$ a measurable nonsingular flow, and $(a^{(1)}$ a
cocycle for this flow, and the functions $(g_t)$ have full support.
We need to prove that the flow $(\psi_t)$ is positive,  and by part (i)
we need to show that for every $w\in \bbw$ 
\begin{equation} \label{e:also.inf}
\int_T w(t) |g_t|^\alpha \lambda(dt) = \infty \ \ \text{$m_1$-a.e.}.
\end{equation}
Let 
$$
{\tt P}_{w}^{f} = \left\{ x\in E:\, \int_T w(t+k) |f|^\alpha\circ \phi_t(x)
\frac{dm\circ \phi_t}{dm}(x) \lambda(dt)=\infty \ \ \text{for all
$k\in \field{Z}$}\right\}. 
$$
Then ${\tt P}_{w}^{f}$ is a flow invariant set (both in the cases
$T=\field{Z}$ and $T=\field{R}$) and, since for every $k$ $w(\cdot
+k)$ differs from a function in $\bbw$ only on a compact interval, it
has a full measure by part (i). Therefore, one can
restrict the integral in (\ref{e:stab.process}) and the flow to 
${\tt P}_{w}^{f}$. By Theorem 1.1 in \cite{rosinski:1995}, there are
measurable maps 
$\Phi: S\to {\tt P}_{w}^{f}$ and $h:S\to \field{R}-\{0\}$ in the
real-valued case 
and $h:S\to \field{C}-\{0\}$ in the complex-valued case such that 
$$
g_t(y) = h(y) f_t(\Phi(y)) \ \ \text{for $m_1\times\lambda$-a.e.
$(y,t)$.}
$$
Choose a measurable set $S_1\subset S$ of full measure such that for
every $y\in S_1$ this relation holds for $\lambda$-a.e. $t$. 
We have for any $y\in S_1$ 
$$
\int_T w(t) |g_t(y)|^\alpha \lambda(dt) = 
|h(y)| \int_T w(t) \left| f_t(\Phi(y))\right|^\alpha \lambda(dt)
= \infty
$$
since $\Phi(y)\in {\tt P}_{w}^{f}$. Therefore, $(\psi_t)$ is positive.

The argument for null flows is similar. 

\end{proof}

Similarly to the conservative-dissipative decomposition, one can check
conditions of the type given in parts (i) and (ii) of Theorem
\ref{t:invariant} without having to have a specific form of integral
representation of a stationary process. Here we are allowing any
{\it measurable} representation (\ref{e:stab.process}). This is a
representation in which the function $f_t(x)$ is jointly measurable in
the variables $x\in E$ and $t\in T$. Every measurable \SaS\ process
has such representation; see Section 11.1 in
\cite{samorodnitsky:taqqu:1994}. In addition, we will assume the full
support condition (\ref{e:full.support}). 
\begin{corollary} \label{c:any.form}
A stationary  \SaS\ with a measurable and full support representation
(\ref{e:stab.process}) is generated by a positive flow 
if and only for any $w\in \bbw$ we have
\begin{equation} \label{e:w.int.infinite.gen}
\int_T w(t) |f_t|^\alpha \lambda(dt)=\infty \ \ \text{$m$-a.e.} 
\end{equation}

A process is generated by a null flow if and
 only for some $w\in \bbw$ we have 
\begin{equation} \label{e:w.int.finite.gen}
\int_T w(t) |f_t|^\alpha \lambda(dt)<\infty \ \ \text{$m$-a.e.} 
\end{equation}
\end{corollary}
\begin{proof}
The same argument as in the proof of Corollary 4.2 in
\cite{rosinski:1995} reduces the situation to that in Theorem
\ref{t:invariant}. 
\end{proof}

\begin{remark} \label{rk:direct.no.full}
{\rm
Note that the direct parts of Corollary \ref{c:any.form} holds without
the assumption of full support. That is, if
(\ref{e:w.int.infinite.gen}) (respectively,
(\ref{e:w.int.finite.gen})) holds then the process is generated by a
positive (respectively, null) flow. Indeed, the above statements will
still be valid if one reduces the integration to the support of the
kernel $(f_t)$. 
}
\end{remark}

Given a stationary \SaS\ process $\BX$ and any integral representation
(\ref{e:stab.process}) of the form (\ref{e:station.kernel}) and of
full support, use the positive-null decomposition $E={\tt P}\cup {\tt
N}$ to write the following analog of (\ref{e:diss-cons}): 
\begin{equation} \label{e:pos-null}
X(t) = \int_{\tt P} f_t(x) \; M(dx) + \int_{\tt N} f_t(x) \; M(dx)
:= X^P(t) + X^N(t), \; t\in T\, 
\end{equation}
a sum of two independent stationary \SaS\ processes, one of which is
generated by a positive flow, and the other one by a null flow.
The next result shows that this decomposition is unique. 
\begin{theorem} \label{t:pos-null}
A decomposition of a stationary \SaS\ process into a sum of two
independent stationary \SaS\ processes, one of which is generated by a
positive flow, and the other one by a null flow, is unique in law.
\end{theorem}
\begin{proof} Consider any fixed {\it minimal} representation of the process
$\BX$ 
$$
X(t) = \int_S g_t(y) \; M_1(dy), \; t\in T,
$$
where $M_1$ has a control measure $m_1$ and 
$$
g_t(y) = a^{(1)}_t(y) \,\left( \frac{dm_1\circ
\psi_t}{dm_1}(y)\right)^{1/\alpha}  g\circ \psi_t(y)
$$
with $(\psi_t)$ a measurable nonsingular flow, and $(a^{(1)}$ a
cocycle for this flow, and the functions $(g_t)$ have full support. We
refer the reader to \cite{rosinski:1995} about basic facts and
properties on minimal representations. Let 
$\BX=\BX_\psi^P+\BX_\psi^N$ be the decomposition in (\ref{e:pos-null})
with respect to that representation.

Consider now any full support representation (\ref{e:stab.process}) of
the form 
(\ref{e:station.kernel}) of the process $\BX$, and let 
$\BX=\BX^P+\BX^N$ be the corresponding decomposition. We will prove
that this decomposition coincides with the above decomposition with
respect to the fixed minimal representation. By Remark 2.5 in
\cite{rosinski:1995} there is a measurable map $\Phi: E\to S$ and a
non-vanishing (real or complex-valued) measurable function $h$ such
that for each $t\in T$, $f_t(x) = h(x) g_t(\Phi(x))$ $m$-a.e., and 
$m \sim m_1\circ \Phi^{-1}$. Let $E={\tt P}\cup {\tt
N}$ be the  positive-null decomposition of the flow $(\phi_t)$ and 
$S={\tt P}^\psi\cup {\tt N}^\psi$ be the  positive-null decomposition
of the flow $(\psi_t)$. Let us prove that
\begin{equation} \label{e:PtoP}
{\tt P} = \Phi^{-1}\left( {\tt P}^\psi\right), \ \ 
{\tt N} = \Phi^{-1}\left( {\tt N}^\psi\right)
\end{equation}
modulo sets of $m$-measure zero.

If $T=\field{Z}$, select (using Theorem 3 in \cite{krengel:1967}) a
$w\in\bbw$ such that 
$$
 {\tt P}^\psi = \left\{ y\in S:\ 
\sum_{n=-\infty}^\infty w(n) |g\circ \psi_n(y)|^\alpha \frac{dm_1\circ
\psi_n}{dm_1}(y) = \infty\right\}
$$
up to a set of $m_1$-measure zero. Notice (once again, by Theorem 3 in
\cite{krengel:1967}) that for $m$-a.e. $x\in {\tt P}$, 
$$
\sum_{n=-\infty}^\infty w(n) |f\circ \phi_n(x)|^\alpha \frac{dm\circ
\phi_n}{dm}(x) = \infty.
$$
Therefore also 
$$
\sum_{n=-\infty}^\infty w(n) |g\circ \psi_n(\Phi(x))|^\alpha \frac{dm_1\circ
\psi_n}{dm_1}(\Phi(x)) = \infty, 
$$
and so $\Phi(x)\in {\tt P}^\psi$ for $m$-a.e. $x\in {\tt P}$. Hence
\begin{equation} \label{e:PtoP1}
\Phi^{-1}\left( {\tt N}^\psi\right) \subset {\tt N} \ \ \text{modulo a
set of $m$-measure zero.}
\end{equation}

Select now a (perhaps different) $w\in\bbw$ such that 
$$
 {\tt P} = \left\{ x\in E:\ 
\sum_{n=-\infty}^\infty w(n) |f\circ \phi_n(x)|^\alpha \frac{dm\circ
\phi_n}{dm}(x) = \infty\right\}
$$
up to a set of $m$-measure zero. Then for $m_1$-a.e. $y\in {\tt
P}^\psi$, 
$$
\sum_{n=-\infty}^\infty w(n) |g\circ \psi_n(y)|^\alpha \frac{dm_1\circ
\psi_n}{dm_1}(y) = \infty.
$$
Since for $m$-a.e. $x\in {\tt N}$ 
$$
\sum_{n=-\infty}^\infty w(n) |g\circ \psi_n(\Phi(x))|^\alpha \frac{dm_1\circ
\psi_n}{dm_1}(\Phi(x)) < \infty, 
$$
we see that $\Phi(x)\in {\tt N}^\psi$ for $m$-a.e. $x\in {\tt N}$. Hence
\begin{equation} \label{e:PtoP2}
\Phi^{-1}\left( {\tt P}^\psi\right) \subset {\tt P} \ \ \text{modulo a
set of $m$-measure zero,}
\end{equation}
and (\ref{e:PtoP}) follows from (\ref{e:PtoP1}) and (\ref{e:PtoP2}). 

In the case $T=\field{R}$, construct the functions $g_+$ and $g_-$ as
in (\ref{e:gplus}) and (\ref{e:gminus}) for the flow $(\phi_t)$ and
the corresponding functions $g_+^\psi$ and $g_-^\psi$ for the flow
$(\psi_t)$. Notice that $g_+(x) = |h(x)|^\alpha g_+^\psi(\Phi(x))$ and
$g_-(x) = |h(x)|^\alpha g_-^\psi(\Phi(x))$ $m$-a.e. This reduced
everything to the discrete time case. For example, there a $w\in\bbw$
such that
$$
 {\tt P}^\psi = \left\{ y\in S:\ 
\sum_{n=0}^\infty w(n) |g_+^\psi\circ \psi_n(y)| \frac{dm_1\circ
\psi_n}{dm_1}(y) + \sum_{n=-\infty}^0 w(n) |g_-^\psi\circ
\psi_n(y)| \frac{dm_1\circ \psi_n}{dm_1}(y)
= \infty\right\}
$$
up to a set of $m_1$-measure zero, etc, and so we obtain
(\ref{e:PtoP}) as in the discrete time case.

Once (\ref{e:PtoP}) has been established, the claim $\BX^P\eid
\BX_\psi^P$ and  
$\BX^N\eid \BX_\psi^N$ follows as in the proof of Theorem 4.3 in
\cite{rosinski:1995}. 
\end{proof}

\begin{remark} \label{rk:comp.decomp}
{\rm 
Notice the relationship between the conservative-dissipative
decomposition (\ref{e:diss-cons}) and the positive-null decomposition
(\ref{t:pos-null}) of a stationary \SaS\ process. Since every
dissipative flow is null, we have a unique in law decomposition
\begin{equation} \label{e:decomp.three}
X(t) = X^D(t) + X^{CN}(t) + X^{P}(t), \ t\in T,
\end{equation} 
into a sum of three independent stationary \SaS\ processes: one
generated by a dissipative flow (the dissipative component of $\BX$),
one generated by a conservative null flow (the conservative null part
of $\BX$) and one generated by a  positive flow (the
positive part of $\BX$). In particular, we have the
relations
\begin{equation} \label{e:decomp.null}
X^N(t) = X^D(t) + X^{CN}(t), \ t\in T,
\end{equation} 
and
\begin{equation} \label{e:decomp.cons}
X^C(t) = X^P(t) + X^{CN}(t), \ t\in T.
\end{equation} 
}
\end{remark}

\begin{remark} \label{rk:structure}
{\rm 
More explicit descriptions of different components in the
decompositions described in Remark \ref{rk:comp.decomp} are desirable.
For the dissipative component $\BX^D$ we have the mixed moving average
representation of \cite{rosinski:1995}: 
\begin{equation} \label{e:mixed.MA}
X^D(t) = \int_W\int_\bbr f(v,x-t)\; M(dv,dx)\,,\quad t\in\reals,
\end{equation}
with $M$ a \SaS\ random measure on a product
measurable space 
$(W\times \bbr,\protect{\bbw}\times\protect{\bbb} )$ 
with control measure $m=\nu\times \text{Leb}$, where $\nu$ is a
$\sigma$-finite measure on $(W, \protect{\bbw})$, and $f\in
L^\alpha(m,\protect{\bbw}\times\protect{\bbb})$. This representation
has proven to be very useful in investigations of properties of
processes generated by dissipative flows; see e.g.
\cite{mikosch:samorodnitsky:2000a}, \cite{resnick:samorodnitsky:2004},
\cite{samorodnitsky:2004a} and \cite{samorodnitsky:2004c}. 

Consider now the positive component $\BX^P$. Assuming this component
is not degenerate (identically zero),  take any representation 
(\ref{e:stab.process}) of the form (\ref{e:station.kernel}), and
recall that there is a probability measure $\mu$ on $\tt P$ invariant
under the flow $(\phi_t)$ and equivalent to the restriction of the
control measure $m$ to $\tt P$. Therefore, we can write, in law,
\begin{equation} \label{e:positive}
X^P(t) = \int_{\tt P} a_t(x) g\circ \phi_t(x) \, M_\mu(dx), \; t\in T\, 
\end{equation}
where $M_\mu$ is a \SaS\ random measure on $E$ with control measure
$\mu$, $(a_t)$ is a length 1 cocycle for the flow $(\phi_t)$, and
$$
g(x) = f(x) \left( \frac{dm}{d\mu}(x)\right)^{1/\alpha}, \ x\in E.
$$
This leads us to the following representation of stationary \SaS\
processes generated by positive flows.

Let $\bigl(\Omega^\prime, {\mathcal F}^\prime, P^\prime\bigl)$ be a
probability space, and $(\theta_t)_{t\in T}$ a measurable group of
shift operators (maps of $\Omega^\prime$ into itself.) Recall that a
stochastic process $\BA = (A(t),\, t\in T)$ is called a (raw)
multiplicative functional if for every $s,t\in \field{R}$ 
$$
A(t) = A(s)A(t)\circ \theta_s \ \ \text{$P^\prime$-a.s.,}
$$
see e.g. \cite{sharpe:1988}. 
Let $\BA$ be  a multiplicative process with $|A(t)|=1$ for all $t\in
T$. 
Let $\BB=(B(t), \, t\in T) = (B(0)\circ \theta_t,\, t\in T) $ be a
compatible with the shift $(\theta_t)$ 
stationary process on $\bigl(\Omega^\prime, {\mathcal F}^\prime,
P^\prime\bigl)$,  such that  
$E^\prime |B(0)|^\alpha<\infty$. 

Let, finally, $M$ be a \SaS\ random measure on  $\bigl(\Omega^\prime,
{\mathcal F}^\prime)$ with control measure $P^\prime$ (defined on a
generic probability space $\bigl(\Omega, {\mathcal F}, P\bigl)$). 
Then 
\begin{equation} \label{e:canonic.pos}
X(t) = \int_{\Omega^\prime} A(t)B(t)\, dM^\prime, \ t\in T
\end{equation}
is, clearly,  a symmetric \SaS\ process (on $\bigl(\Omega, {\mathcal F},
P\bigl)$) corresponding to a positive flow. Furthermore, the
expression (\ref{e:positive}) shows that every process corresponding 
to a positive flow has a representation of the form
(\ref{e:canonic.pos}). 

The representation in (\ref{e:canonic.pos}) describes the positive
component of a stationary \SaS\ process in terms of an stationary
process with a finite moment of order $\alpha$ and a multiplicative
functional of absolute value 1. It is an attractive representation
from purely probabilistic point of view because both stationary
processes and multiplicative functionals are important and widely
studied objects of their own. 

Several particular cases have been long considered in literature. If
the multiplicative functional $\BA$ is identically equal to 1, then
the \SaS\ process has the form 
\begin{equation} \label{e:doubly.stat}
X(t) = \int_{\Omega^\prime} B(t)\, dM^\prime, \ t\in T,
\end{equation}
an integral of an $L^\alpha$ stationary process. Processes of the form
(\ref{e:doubly.stat}) have been known as {\it doubly stationary
processes} (the term introduced, apparently, by
\cite{cambanis:hardin:weron:1987}) with a finite control measure. The
further particular cases where the stationary process $\BB$ is a
stationary zero mean Gaussian process or a symmetric $\beta$-stable
process with $\alpha<\beta<2$ are known as, correspondingly,
sub-Gaussian or sub-stable processes and have an alternative
representation as
\begin{equation} \label{e:sub-stab}
X(t) = c_{\alpha,\beta} W^{1/\beta} B(t),\ t\in T,
\end{equation}
where now $\BX$ and $\BB$ live on the same probability space, $
c_{\alpha,\beta}$ is a positive constant for $\alpha<\beta\leq 2$, and
$W$ is an independent of $\BB$ positive strictly $\alpha/\beta$-stable
random variable. See \cite{samorodnitsky:taqqu:1994}. 

Another important example of processes generated by positive flows
that has been long studied is that of {\it harmonizable processes}. It
corresponds to taking the stationary
process $\BB$ to be constant and the multiplicative functional $\BA$ 
of a special form. See \cite{rosinski:1995}.

No ``canonical'' representation of the component $\BX^{CN}$ seems to
be known at this time!
}
\end{remark}

\section{Ergodicity} \label{s:ergodicity}

Let $\bigl(\Omega, {\mathcal F}, P\bigl)$) be a probability space, and
$(\theta_t)_{t\in T}$ a measurable group of shift operators preserving
$P$. Let $X(t)=X(0)\circ \theta_t, \, t\in\field{R}$ be a stationary
process. Recall that a process $\BX$ is {\it ergodic} if for every events
$A,B\in \sigma(X(t), \, t\in\field{R})$ 
$$
\lim_{T\to\infty}\frac1T \int_0^T P\left( A\cap B\circ \theta_t\right)\,
\lambda(dt) = P(A)P(B).
$$
A process $\BX$ is {\it weakly mixing} if for all $A$ and $B$ as above
$$
\lim_{T\to\infty}\frac1T \int_0^T \left| P\left( A\cap B\circ
\theta_t\right) - P(A)P(B)\right| \, \lambda(dt) = 0.
$$
In general, weak mixing is a stronger condition than ergodicity, but
the two notions coincide for both stationary Gaussian processes (see
\cite{maruyama:1949}) and for stationary infinitely divisible
processes, see \cite{rosinski:zak:1997} (for stationary 
\SaS\ processes  this statement was established by 
\cite{podgorski:1992}). Furthermore, a stationary Gaussian processes
is ergodic (weakly mixing) if and only if its spectral measure is atomless
(again, \cite{maruyama:1949}). 

The following result shows that for stationary \SaS\ processes with
$0<\alpha<2$ ergodicity is related to absence of the  component
corresponding to a positive flow.
\begin{theorem} \label{t:ergodic}
A stationary \SaS\ process is ergodic (equivalently, weakly mixing) if
and only if the component 
$\BX^P$ in (\ref{e:decomp.three}) corresponding to a positive flow
vanishes. 
\end{theorem}
\begin{proof}
We start with proving that a process generated by a null flow is
weakly mixing. Consider first the case $T=\field{Z}$. Let $\BX$ be
given in the form (\ref{e:stab.process}) with the kernel of the form
(\ref{e:station.kernel}). By
\cite{gross:1994} (see also the discussion of weak mixing in
\cite{krengel:1967}) weak mixing is equivalent to the following
statement: for every compact set $K$ bounded
away from zero and $\epsilon>0$ 
\begin{equation} \label{e:weak.mix1}
\lim_{n\to\infty} \frac1n \sum_{j=0}^{n-1} m\left( x:\,
|f_0(x)|^\alpha \in K,\, |f_j(x)|^\alpha>\epsilon\right) = 0.
\end{equation}
Fix $K$ as above and let $A= \{ x:\, |f_0(x)|^\alpha \in K\}$. Then
$m(A)<\infty$. For $\theta>0$ we have
$$
\frac1n \sum_{j=0}^{n-1} m\left( x\in A:\,
|f_j(x)|^\alpha>\epsilon\right)
= \int_A \frac1n \sum_{j=0}^{n-1} \one\left( x\in A:\,
|f_j(x)|^\alpha>\epsilon\right)\, m(dx)
$$
$$
= \int_A \left(\frac1n \sum_{j=0}^{n-1} \one\left( x\in A:\,
|f_j(x)|^\alpha>\epsilon\right)\right) \one\left( \frac1n
\sum_{j=0}^{n-1} |f_j(x)|^\alpha>\theta\right)\, m(dx)
$$
\begin{equation} \label{e:bound.null}
+ \int_A \left(\frac1n \sum_{j=0}^{n-1} \one\left( x\in A:\,
|f_j(x)|^\alpha>\epsilon\right)\right) \one\left( \frac1n
\sum_{j=0}^{n-1} |f_j(x)|^\alpha\leq \theta\right)\, m(dx)
\end{equation}
$$
\leq m\left( x\in A:\, \frac1n
\sum_{j=0}^{n-1} |f_j(x)|^\alpha>\theta\right)
+ \epsilon^{-1} \int_A \frac1n
\sum_{j=0}^{n-1} |f_j(x)|^\alpha
\one\left( \frac1n
\sum_{j=0}^{n-1} |f_j(x)|^\alpha\leq \theta\right)\, m(dx)
$$
$$
\leq m\left( x\in A:\, \frac1n
\sum_{j=0}^{n-1} |f_j(x)|^\alpha>\theta\right)
+ \frac{\theta}{\epsilon}m(A).
$$
Let $T:\, L^1(m)\to  L^1(m)$ be defined by 
$$
Tg =  \frac{dm\circ \phi}{dm} g\circ \phi_1\,.
$$
Note that $T$ is a positive contraction (actually, an isometry) on
$L^1(m)$. Furthermore, $|f_j|^\alpha=T^j|f|^\alpha$ for all $j\geq 0$. 

By the stochastic ergodic theorem (see Theorem 4.9, page 143 in
\cite{krengel:1985}) the average 
$$
\frac1n
\sum_{j=0}^{n-1} |f_j|^\alpha =
\frac1n
\sum_{j=0}^{n-1} T^j|f|^\alpha
$$
converges in measure, on any set of a finite measure, to a nonnegative
function $h\in L^1(m)$ which invariant under $T$. This means that the finite
measure $d\mu = h\, dm$ is invariant under the map $\phi$. Since the
flow is null, the measure $\mu$ must be the null measure and, hence,
the limit function $h$ must be equal to zero. The convergence in
measure then gives us that 
$$
\lim_{n\to\infty} m\left( x\in A:\, \frac1n
\sum_{j=0}^{n-1} |f_j(x)|^\alpha>\theta\right) = 0
$$
for every $\theta>0$. Letting in (\ref{e:bound.null}) first
$n\to\infty$,   and then letting  $\theta\to 
0$ verifies (\ref{e:weak.mix1}) and, hence, shows that a process
generated by a null flow is weakly mixing. 

Suppose now that $T=\field{R}$. By Theorem 2 of \cite{rosinski:zak:1997}
ergodicity of $\BX$ is equivalent to 
\begin{equation} \label{e:weak.mix2}
\lim_{M\to\infty} \frac1M \int_0^M \exp\left\{ 2\| X(0)\|^\alpha_\alpha -
\| X(t)-X(0)\|^\alpha_\alpha\right\}\, dt =1. 
\end{equation}
Here $\| X(t)\|_\alpha$ is simply the scaling parameter of the 
\SaS\ random variable $X(t)$; see \cite{samorodnitsky:taqqu:1994}. 

Recall that our processes are continuous in probability. Therefore, 
for every $\epsilon>0$ we can choose $\delta>0$ such that $\|
X(t)-X(s)\|_\alpha\leq \epsilon$ if $|t-s|\leq \delta$. For such a
pair of $\epsilon$ and $\delta$ we have for $N=1,2,\ldots $
\begin{equation} \label{e:weak.mix3}
\frac{1}{N\delta}\int_0^{N\delta} \exp\left\{ 2\| X(0)\|^\alpha_\alpha -
\| X(t)-X(0)\|^\alpha_\alpha \right\}\, dt 
= \frac{1}{N\delta} \sum_{j=1}^N I_j(\delta),
\end{equation}
where 
$$
I_j(\delta) = \int_{(j-1)\delta}^{j\delta}
\exp\left\{ 2\| X(0)\|^\alpha_\alpha -
\| X(t)-X(0)\|^\alpha_\alpha\right\}\, dt,
$$
$j=1,\ldots, N$. 

Suppose that $1< \alpha<2$. Then $\|\cdot
\|_\alpha$ satisfies the triangle inequality, and by the inequality
$$
(a+b)^\alpha \leq a^\alpha + \alpha ab^{\alpha-1} +b^\alpha
$$
for $a,b\geq 0$ we have for every $j=1,\ldots, N$ and $t\in
[(j-1)\delta, j\delta]$, 
$$
\| X(t)-X(0)\|^\alpha_\alpha \leq \left( \epsilon+ \|
X\bigl((j-1)\delta\bigr) -X(0)\|_\alpha\right)^\alpha 
$$
$$
\leq \| X\bigl((j-1)\delta\bigr) -X(0)\|_\alpha^\alpha +\alpha\epsilon
\| X\bigl((j-1)\delta\bigr) -X(0)\|_\alpha^{\alpha-1} +\epsilon^\alpha.
$$
Therefore,
$$
I_j(\delta) \geq \delta \exp\left\{ 2\| X(0)\|^\alpha_\alpha -
\| X\bigl((j-1)\delta\bigr) -X(0) \|_\alpha^\alpha\right\}
\exp\left\{ -\epsilon^\alpha - \alpha\epsilon
\| X\bigl((j-1)\delta\bigr) -X(0)\|_\alpha^{\alpha-1}\right\}.
$$
Therefore, 
$$
\liminf_{M\to\infty} \frac1M \int_0^M \exp\left\{ 2\| X(0)\|^\alpha_\alpha -
\| X(t)-X(0)\|_\alpha^\alpha\right\}\, dt
$$
$$
\geq \exp\left\{ -\epsilon^\alpha - 2\alpha\epsilon
\|X(0)\|_\alpha^{\alpha-1}\right\}
\liminf_{N\to\infty} \frac1N \sum_{j=1}^N
\exp\left\{ 2\| X(0)\|^\alpha_\alpha -
\| X\bigl((j-1)\delta\bigr) -X(0) \|_\alpha^\alpha\right\}.
$$
Since the flow $(\phi_t)$ is null, so is the discrete time flow
$(\phi_n)$ and, hence, the discrete time process $\bigl( X(n), \, n\in
\field{Z}\bigr)$ is ergodic. Therefore, 
$$
\lim_{N\to\infty} \frac1N \sum_{j=1}^N
\exp\left\{ 2\| X(0)\|^\alpha_\alpha -
\| X\bigl((j-1)\delta\bigr) -X(0) \|_\alpha^\alpha\right\} =1.
$$
Letting $\epsilon\to 0$ we conclude that 
$$
\liminf_{M\to\infty} \frac1M \int_0^M \exp\left\{ 2\| X(0)\|^\alpha_\alpha -
\| X(t)-X(0)\|_\alpha^\alpha \right\}\, dt \geq 1.
$$
Similarly, using the inequality
$$
(a-b)_+^\alpha \geq a^\alpha-\alpha ab^{\alpha -1}
$$
for $1<\alpha<2$ and $a,b\geq 0$ we have
for every $j=1,\ldots, N$ and $t\in
[(j-1)\delta, j\delta]$, 
$$
\| X(t)-X(0)\|^\alpha_\alpha \geq \left( \| X\bigl((j-1)\delta\bigr)
-X(0) \|_\alpha-\epsilon\right)_+^\alpha  
$$
$$
\geq \| X\bigl((j-1)\delta\bigr) -X(0) \|_\alpha^\alpha
-\alpha\epsilon^{\alpha-1}\| X\bigl((j-1)\delta\bigr)
-X(0) \|_\alpha,
$$
and the same argument as above gives us 
$$
\limsup_{M\to\infty} \frac1M \int_0^M \exp\left\{ 2\|
X(0)\|^\alpha_\alpha - \| X(t)-X(0)\|^\alpha_\alpha\right\}\, dt \leq 1.
$$
Therefore, (\ref{e:weak.mix2}) holds. The argument for
(\ref{e:weak.mix2}) in the case $0<\alpha\leq 1$ is similar (and even
easier). Hence, a process
generated by a null flow is ergodic (weakly mixing) in the case
$T=\field{R}$ as well. 

In the opposite direction, suppose that the component $\BX^P$ in
(\ref{e:decomp.three}) corresponding to a positive flow does not 
vanish. The fact that $\BX^P$ is not ergodic follows from
Corollary 4.2 in \cite{gross:1994}. Therefore, there are two
stationary processes with different finite-dimensional distributions,
$\BY^{(1)}$ and $\BY^{(2)}$ and $0<p<1$ , such that 
$$
\BX^P = \left\{
\begin{array}{ll}
\BY^{(1)} & \text{with probability $p$}   \\[2mm]
\BY^{(2)} & \text{with probability $1-p$} 
\end{array}
\right. 
$$
(to see this, view the process $\BX^P$ as defined on the canonical
path space, take an invariant event $A$ with probability $0<p<1$, and
let $\BY^{(1)}$ have the law of $\BX^P$ conditioned on belonging to
$A$, and let $\BY^{(2)}$ have the law of $\BX^P$ conditioned on belonging to
$A^c$.) 

Therefore, by (\ref{e:pos-null}),
\begin{equation} \label{e:weak.mix.last}
\BX = \left\{
\begin{array}{ll}
 \BX^N+\BY^{(1)} & \text{with probability $p$}   \\[2mm]
 \BX^N+\BY^{(2)} & \text{with probability $1-p$} 
\end{array}
\right. .
\end{equation}
Since the characteristic functions of infinitely divisible random
vectors (and, in particular, of \SaS\ random vectors) do not vanish,
we see that the processes $ \BX^N+\BY^{(1)}$ and $ \BX^N+\BY^{(2)}$
have different finite dimensional distributions. Therefore,
(\ref{e:weak.mix.last}) implies that $\BX$ is not ergodic. 
\end{proof}
\begin{remark} \label{rk:erg.not.mix}
{\rm
It is known (see \cite{surgailis:rosinski:mandrekar:cambanis:1993})
that all stationary \SaS\ processes generated by dissipative flows are
mixing, and we have seen that all processes generated by positive
flows are not ergodic. On the other hand, there are processes
generated by conservative null flows that are ergodic but not mixing,
simply because ergodicity and mixing are not equivalent for stationary
\SaS\ processes; see \cite{gross:robertson:1993}. 
}
\end{remark}

Combining the results of Theorem \ref{t:ergodic} above and of Theorem
4.1 in 
\cite{samorodnitsky:2004a} we obtain the following interesting
characterization of processes corresponding to  conservative null
flows. 
\begin{corollary} \label{c:max.erg}
A stationary \SaS\ process $(X_n,\, n\in \field{Z})$ is generated by a
conservative null flow if and only if it is ergodic and
$n^{-1/\alpha}\max_{j=1,\ldots ,n}|X_j|\to 0$ is probability as
$n\to\infty$. 
\end{corollary}

A similar characterization holds for locally bounded continuous time
stationary \SaS\ processes (see \cite{samorodnitsky:2004c}).

\section{Examples} \label{s:examples}

In this section we consider several examples of both discrete and
continuous time stationary \SaS\ processes and see how the notions
developed in the previous sections apply here. 

\begin{example} \label{ex:null-rec.MC}
{\rm
Consider a bilateral real-valued Markov chain 
with law $P_x(\cdot)$, $x\in\field{R}$, admitting an infinite
$\sigma$-finite invariant measure $\pi$.  Define a
$\sigma$-finite measure $m$ on $E=\field{R}^{\field{Z}}$ by 
$$
m(\cdot)= \int_{\field{R}}
P_x(\cdot)\, \pi(dx)\,. 
$$
Suppose that the Markov chain is $m$-irreducible and Harris recurrent
(see e.g. \cite{meyn:tweedie:1993}. 

Let $\BX$ be a stationary discrete-time \SaS\ process defined by the
integral  
representation (\ref{e:stab.process}), with $M$ being a \SaS\ random
measure with control measure $m$, and 
$$
f_n(x)= a_n(x)\,f\circ \phi_n(x)\, \ x\in E\,, \ n=0,1,2,\ldots\,,
$$
where $\phi$ is the left shift operator on $E$, $f\in L^{\alpha}(m)$,
and $a$ a $\pm 1$-valued cocycle for the flow $(\phi_n)$. 

This process was considered in \cite{rosinski:samorodnitsky:1996} who 
showed that this process is generated by a conservative flow and is
mixing. By Theorem \ref{t:ergodic} we conclude that this process
corresponds to a conservative null flow. 
}
\end{example}

\begin{example} \label{ex:cycle}
{\rm
Here we consider the class of continuous time stationary \SaS\
processes corresponding to the so-called {\it cyclic flows} introduced
by \cite{pipiras:taqqu:2004}. These processes have a representation of
the form
\begin{equation} \label{e:cyclic.rep}
X(t) = \int_S\int_{[0,q(z))} b(z)^{[v+s(z)t]_{q(z)}}\, g\left(
z,\left\{ v+s(z)t\right\}_{q(z)}\right)\, M(dz,dv), \ t\in \field{R},
\end{equation}
where $\bigl( S, {\mathcal S}\bigr)$ is a measurable space, $b$, $q$
and $s$ are measurable (possibly, complex-valued) functions on $S$,
such that $|b(z)|=1$, $q(z)>0$ and $s(z)\not= 0$ for all $z\in S$.
Furthermore, $M$ is a \SaS\, random measure on $E=\bigl\{ (z,v)\in
S\times \field{R}_+:\, 0\leq v<q(z)\bigr\}$ with control measure
$m=\sigma\times\lambda$, where $\sigma$ is a $\sigma$-finite measure on
$S$. Finally, $g\in L^\alpha(m)$ is assumed to be such that 
\begin{equation} \label{e:cyclic.full}
\sigma\left\{ z\in S:\, g(z,\cdot) = 0 \ \ \text{a.e. on
$[0,q(z))$}\right\} = 0
\end{equation}
(the full support assumption). In (\ref{e:cyclic.rep}), $[\cdot]_a$
and $\{ \cdot\}_a$ are 
respectively, the integer part and the fractional part with respect to
a positive number $a$. 

We will see that such a process is generated by a positive flow. We
use Corollary \ref{c:any.form}. Assume e.g. that $s(z)>0$. For any
$w\in \bbw$  
$$
\int_{-\infty}^\infty w(t) |f_t(z,v)|^\alpha \lambda(dt)
= \int_{-\infty}^\infty w(t) \left| g\left(
z,\left\{ v+s(z)t\right\}_{q(z)}\right)\right|^\alpha \lambda(dt)
$$
\begin{equation} \label{e:cyclic.3}
= \frac{1}{s(z)} \int_{-\infty}^\infty w\left( \frac{t-v}{s(z)}\right)
\left| g\left(
z,\left\{t\right\}_{q(z)}\right)\right|^\alpha \lambda(dt)
\end{equation}
$$
= \frac{1}{s(z)} \sum_{n=-\infty}^\infty \int_0^{q(z)} 
w\left( \frac{t-v+nq(z)}{s(z)}\right)\left| g(z,t)\right|^\alpha \lambda(dt)
$$
$$
= \frac{1}{s(z)}\int_0^{q(z)} \left| g(z,t)\right|^\alpha
\left( \sum_{n=-\infty}^\infty w\left(
\frac{t-v+nq(z)}{s(z)}\right)\right)\lambda(dt).
$$
Notice that for any $v,z$ and $t$
$$
 \sum_{n=-\infty}^\infty w\left(
\frac{t-v+nq(z)}{s(z)}\right)
\geq \int_{(v-t)_+/q(z)+1}^\infty w\left(
\frac{t-v+xq(z)}{s(z)}\right)dx=\infty.
$$
Therefore, by (\ref{e:cyclic.3}) and the full support assumption
(\ref{e:cyclic.full}) we conclude that 
$$
\int_{-\infty}^\infty w(t) |f_t(z,v)|^\alpha \lambda(dt) = \infty
$$
for $m$-a.e. $(z,v)\in E$ such that $s(z)>0$. Since the argument in
the case $s(z)<0$ is the same, by Corollary \ref{c:any.form} the
cyclic process
is generated by a positive flow. 
}
\end{example}

\begin{example} \label{ex:stat.incr}
{\rm
Let $E=\Omega_1\times\field{R}$, where 
$\bigl(\Omega_1, {\mathcal F}_1, P_1\bigl)$ is the canonical
probability space $\Omega_1=\field{R}^{\field{R}}$ with the
cylindrical $\sigma$-field, and $P_1$ such that 
$\bigl( Y(t,\omega_1) = \omega_1(t),\, t\in\field{R},\,
\omega_1\in \Omega_1\bigr)$ is a measurable process with stationary
increments.

Let $M$ be a \SaS\ random measure on $E$ with control measure
$P_1\times\lambda$ (where $\lambda$ is, as usual, the Lebesgue measure
on $\field{R}$). Let $\varphi\in L^\alpha(\lambda)$ and
\begin{equation} \label{e:stat.incr.rep}
X(t) = \int_{\Omega_1} \int_{-\infty}^\infty \varphi\left(
Y(t,\omega_1)+z\right)\, M(d\omega_1, dz), \ t\in \field{R}.
\end{equation}
Notice that
$$
\int_{\Omega_1} \int_{-\infty}^\infty \left| \varphi\left(
Y(t,\omega_1)+z\right)\right|^\alpha\, m(d\omega_1, dz)
$$
$$
= \int_{-\infty}^\infty E_1\left| \varphi\left(
Y(t)+z\right)\right|^\alpha\, \lambda(dz) 
= \int_{-\infty}^\infty \varphi(z)\, \lambda(dz)<\infty
$$
for every $t\in \field{R}$, and so the process $\BX$ in
(\ref{e:stat.incr.rep}) is well defined. Furthermore, let 
$(\theta_t)_{t\in \field{R} }$ be the measurable group of left shift
operators on $\Omega_1$, and let $\phi_t(\omega_1,z) =
(\theta_t(\omega_1),z)$ be 
the induced flow on $E$. Notice that the measure $m$ is preserved
under this flow. To see this notice that for every $t_1<t_2<\ldots
<t_k$ and $t>0$ 
$$
m\left\{ (\omega_1,z)\in E:\, Y(t_1+t,\omega_1)+z\in A,\, 
Y(t_2+t,\omega_1)-Y(t_1+t,\omega_1)\in B_1,\ldots
, \right.
$$
$$
\left. 
Y(t_k+t,\omega_1)-Y(t_{1}+t,\omega_1)\in B_{k-1}\right\}
$$
$$
= E_1\left[\one\left( Y(t_2+t,\omega_1)-Y(t_1+t,\omega_1)\in B_1,\ldots
,Y(t_k+t,\omega_1)-Y(t_{1}+t,\omega_1)\in B_{k-1}\right)
\right.
$$
$$
\left.
\int_{-\infty}^\infty \one\left( Y(t_1+t,\omega_1)+z\in A\right)\,
\lambda(dz) \right]
$$
$$
= \lambda(A) P_1\left( Y(t_2+t,\omega_1)-Y(t_1+t,\omega_1)\in B_1,\ldots
,Y(t_k+t,\omega_1)-Y(t_{1}+t,\omega_1)\in B_{k-1}\right)
$$
$$
= \lambda(A) P_1\left( Y(t_2,\omega_1)-Y(t_1,\omega_1)\in B_1,\ldots
,Y(t_k,\omega_1)-Y(t_{1},\omega_1)\in B_{k-1}\right)
$$
$$
= m\left\{ (\omega_1,z)\in E:\, Y(t_1,\omega_1)+z\in A,\, 
Y(t_2,\omega_1)-Y(t_1,\omega_1)\in B_1,\ldots
, \right.
$$
$$
\left. 
Y(t_k,\omega_1)-Y(t_{1},\omega_1)\in B_{k-1}\right\}
$$
by the stationarity of the increments of $\BY$,  
for all Borel sets $A$, $B_1,\ldots, B_{k-1}$ with
$\lambda(A)<\infty$. Hence the flow $(\phi_t)$ preserves the measure
of every measurable finite dimensional rectangle, hence it preserves
the measure $m$. Therefore, the process in (\ref{e:stat.incr.rep}) is
a stationary \SaS\ process. 

We will prove that, under the assumption
\begin{equation} \label{e:infty.prob}
|Y(t)|\to \infty \ \ \text{in $P_1$-probability as $|t|\to\infty$}
\end{equation}
the process in (\ref{e:stat.incr.rep}) corresponds to a null flow.

To this end we use Corollary \ref{c:any.form}. It is enough to show
that for some $w\in \bbw$
\begin{equation} \label{e:finite.exp}
\int_{\Omega_1} \int_{-\infty}^\infty \left( \int_{-\infty}^\infty
w(t) \left| \varphi\left(
Y(t,\omega_1)+z\right)\right|^\alpha\, dt\right) \frac{1}{\sqrt{2\pi}}
e^{-z^2/2}\, dz\, P_1(d\omega_1)<\infty.
\end{equation}
Indeed, the integral in (\ref{e:finite.exp}) can be written in the
form $\int_{-\infty}^\infty w(t) h(t)\, dt$, where 
$$
h(t) = E_1\left| \varphi\left(
Y(t)+G\right)\right|^\alpha, \ \ t\in \field{R},
$$
and $G$ is an independent (under $P_1$) of $Y(t)$ standard normal
random variable. Therefore, if we show that
\begin{equation} \label{e:h.to.0}
\lim_{|t|\to\infty} h(t) =0,
\end{equation}
then it is clear how to choose a $w\in \bbw$ such that
(\ref{e:finite.exp}) holds. In order to prove (\ref{e:h.to.0}), note
that for every $M>0$
$$
h(t) = E_1\left[ \left| \varphi\left(
Y(t)+G\right)\right|^\alpha \one\left( |Y(t)|\leq M\right)\right]
+ E_1\left[ \left| \varphi\left(
Y(t)+G\right)\right|^\alpha \one\left( |Y(t)|> M\right)\right]
$$
$$
:= h_1(t;M) + h_2(t;M).
$$
By (\ref{e:infty.prob}) we immediately see that for every $M$, 
$$
h_1(t;M) \leq \frac{1}{\sqrt{2\pi}} \|\varphi\|_\alpha^\alpha
P_1\left(  |Y(t)|\leq M\right) \to 0
$$
as $|t|\to\infty$. On the other hand, by the dominated convergence
theorem, 
$$
\lim_{|y|\to\infty} \int_{-\infty}^\infty \frac{1}{\sqrt{2\pi}}
e^{-(z-y)^2/2} |\varphi(z)|^\alpha\, dz =0.
$$
Therefore, 
$$
\lim_{M\to\infty} \limsup_{|t|\to\infty} h_2(t;M) = 0,
$$
and so (\ref{e:h.to.0}) follows. 

Notice that if (\ref{e:infty.prob}) fails, then the process $\BX$ will
not, in general, be generated by a null flow. For example, if the
process $\BY$ is, actually, stationary under $P_1$, then the process
$\BX$ is generated by a positive flow. To see this notice that, in
this case, the process $\BX$ can be represented in the form
(\ref{e:canonic.pos}) with $\Omega^\prime = \Omega_1\times\field{R}$, 
$P^\prime = P_1\times P_G$, $\BA\equiv 1$ and $B(t)=f(z)^{-1/\alpha}
\varphi(Y(t,\omega_1)+z), \, t\in\field{R}$, where $P_G$ is the
standard Gaussian law on $\field{R}$ and $f$ is the density of this
law. 

Under the assumption (\ref{e:infty.prob}) both conservative null and
dissipative flows are possible. For example, if we strengthen the 
assumption (\ref{e:infty.prob}) to 
\begin{equation} \label{e:infty.a.s.}
|Y(t)|\to \infty \ \ \text{$P_1$-a.s. as $|t|\to\infty$}
\end{equation}
then the process $\BX$ is generated by a dissipative flow. To see this
note that under the assumption (\ref{e:infty.a.s.}) there is a bounded
strictly positive and monotone on both $(-\infty,0]$ and $[0,\infty)$
deterministic function $\tilde\varphi\in 
L^\alpha(\lambda)$ such that 
$$
\int_{-\infty}^\infty \tilde\varphi\left( \frac12 Y(t)\right)\,
dt<\infty
$$
$P_1$-a.s. If we replace the function $\varphi$ in the definition of the
process $\BX$ by $\tilde\varphi$, we immediately conclude by Corollary
4.2 of \cite{rosinski:1995} that the new \SaS\ process is generated by
a dissipative flow. Since the kernel of the new process has support at
least as large as that of the kernel of the process $\BX$, the latter
process is also generated by a dissipative flow. 

On the other hand, if, for example, the process $\BY$ is, under $P_1$, a
fractional Brownian motion with any $0<H<1$, then the process $\BX$ is
generated by a conservative flow (see Section 3 in
\cite{samorodnitsky:2004c}) and, hence, it corresponds
to a conservative null flow. 
}
\end{example}

\bibliographystyle{/net/flow/export/home/gennady/texfiles/mystyle}
\bibliography{/net/flow/export/home/gennady/texfiles/bibfile}

\end{document}